\documentstyle[12pt]{article}

\input amssym.def
\input amssym

\newcommand{\bR}{\bf R}
\newcommand{\bP}{\bf P}
\newcommand{\bB}{\bf B}

\newcommand{\bx}{\bf x}
\newcommand{\bxu}{\underline{\bf x}}
\newcommand{\bxo}{\overline{\bf x}}
\newcommand{\by}{\bf y}
\newcommand{\byu}{\underline{\bf y}}
\newcommand{\byo}{\overline{\bf y}}
\newcommand{\0}{\bf 0}
\newcommand{\1}{\bf 1}

\def\cle{\preccurlyeq}

\def\sumlim{\sum\limits}
\def\maxlim{\max\limits}
\def\suplim{\sup\limits}
\def\CB{\cal B}
\def\CF{\cal F}
\newcommand{\Rmax}{{\bR}_{\max}}
\newcommand{\Rmin}{{\bR}_{\min}}
\newcommand{\Rpl}{{\bR}_+}

\begin{document}
\begin{center}
The 25th Linz Seminar on Fuzzy Set Theory

Linz, Austria, 2004, February 3 -- February 7

\hrule\smallskip
\hrule
\vskip 1cm

{\Large {\bf Dequantization of Mathematics, idempotent 
semirings and fuzzy sets}}

\medskip

{\sc G.~L.~Litvinov}\footnote{Independent University of Moscow,
E-mail: glitvinov@mail.ru

This research was supported by the Russian Foundation for Basic
Research, grant 02--01--01062}
\end{center}

\medskip

{\bf 1. Introduction.} 
  The traditional mathematics over numerical fields can be
dequantized as the Planck constant $\hbar$ tends to zero taking
pure imaginary values.  This dequantization leads to the so-called
Idempotent Mathematics based on replacing the usual arithmetic
operations by a new set of basic operations (e.g., such as maximum
or minimum), that is on the concepts of idempotent semifield and
semiring. Typical examples 
are given by the so-called $(\max, +)$ algebra $\Rmax$ and $(\min, +)$
algebra $\Rmin$. Let $\bR$ be the field of real numbers.
Then ${\Rmax}={\bR} \cup \{-\infty\}$ with operations
$x\oplus y=\max \{x,y\}$ and $x\odot y=x+y$. Similarly
${\Rmin}={\bR}\cup \{+\infty\}$ with the operations 
$\oplus=\min$, $\odot=+$. The new addition $\oplus$ is idempotent,
i.e., $x\oplus x=x$ for all elements $x$. Some problems
that are nonlinear in the traditional sense turn out to be linear
over a suitable idempotent semiring (idempotent superposition 
principle \cite{1}). For example, the Hamilton-Jacobi equation (which is an
idempotent version of the Schr\"odinger equation) is linear over 
$\Rmin$ and~$\Rmax$.

The basic paradigm is expressed in terms of an {\it idempotent 
correspondence principle} \cite{2}.This principle is similar to the 
well-known correspondence
principle of N.~Bohr in quantum theory (and closely related to it).
Actually, there exists a heuristic correspondence between
important, interesting and useful constructions and results of the
traditional mathematics over fields and analogous constructions and results
over idempotent semirings and semifields (i.e., semirings and semifields
with idempotent addition). For example, the well-known Legendre
transform can be treated as an $\Rmax$-version of the traditional
Fourier transform (this observation is due to V.~P.~Maslov).

A systematic and consistent application of the idempotent
correspondence principle leads to a variety of results, often quite
unexpected. As a result, in parallel with the traditional mathematics over
rings, its ``shadow'', the Idempotent Mathematics, appears. This ``shadow''
stands approximately in the same relation to the traditional mathematics as
classical physics to quantum theory. In many respects Idempotent
Mathematics is simpler than the traditional one. However, the transition
from traditional concepts and results to their idempotent analogs is often
nontrivial.

In this talk a brief survey of basic ideas of Idempotent Mathematics is
presented. Relations between this theory and the theory of fuzzy sets
as well as the possibility theory and some applications (including
computer applications) are discussed. Hystorical surveys and the 
corresponding references can be found in \cite{2}--\cite{5}.

{\bf 2. Semirings, semifields, and idempotent dequantization.}
Consider a set $S$ equipped with two algebraic operations: {\it addition} 
$\oplus$ and {\it multiplication} $\odot$. It is a {\it semiring} if the 
following conditions are satisfied:
\begin{itemize}
\item the addition $\oplus$ and the multiplication $\odot$ are
associative;
\item the addition $\oplus$ is commutative;
\item the multiplication $\odot$ is distributive with respect to 
the addition $\oplus$: $x\odot(y\oplus z)=(x\odot y)\oplus(x\odot z)$ and
$(x\oplus y)\odot z=(x\odot z)\oplus(y\odot z)$ for all $x,y,z\in S$.
\end{itemize}

The semiring is {\it commutative} if the multiplication $\odot$ is 
commutative. A {\it unity} of a semiring $S$ is an element ${\1}\in S$ such 
that ${\1}\odot x=x\odot{\1}=x$ for all $x\in S$. A {\it zero} of a semiring
$S$ is an element ${\0}\in S$ such that ${\0}\neq{\1}$ and ${\0}\oplus x=x$, 
${\0}\odot x=x\odot {\0}={\0}$ for all $x\in S$. A semiring $S$ is called an
{\it idempotent semiring} if $x\oplus x=x$ for all $x\in S$. A 
semiring $S$ with neutral elements $\0$ and $\1$ is called a {\it 
semifield} if every nonzero element of $S$ is invertible.

The following examples are important. Let $\bP$ be the segment $[0,1]$
equipped with the operations $x\oplus y = \max \{x, y\}$ and
$x\odot y = \min \{x, y\}$; then $\bP$ is a commutative idempotent
semiring (but not a semifield). The subset ${\bB} =\{0,1\}$ in $M$
equipped with the same operations is the well-known Boolean algebra
which is an idempotent semifield. In this case $\oplus$ and $\odot$
are the usual Boolean operations (disjunction and conjunction).
In the general case the semiring addition and multiplication could be
treated as generalized logical (Boolean) operations.

Let $\bR$ be the field of real numbers and $\Rpl$ the semiring of all
nonnegative real numbers (with respect to the usual addition and 
multiplication). The change of variables $x \mapsto u = h \ln x$, $h > 0$, 
defines a map $\Phi_h \colon \Rpl \to S = {\bR} \cup \{-\infty\}$.  Let the 
addition and multiplication operations be mapped from $\bR$ to $S$ by $\Phi_h$,
i.e., let $u \oplus_h v = h \ln (\exp (u/h) + \exp(v/h))$, $u \odot v = u +
v$, ${\0} = -\infty = \Phi_h(0)$, ${\1} = 0 = \Phi_h(1)$.  It can easily be
checked that $u \oplus_h v \to \max \{u,v\}$ as $h \to 0$ and $S$ forms a
semiring with respect to addition $u \oplus v = \max \{u,v\}$ and
multiplication $u \odot v = u + v$ with zero ${\0} = -\infty$ and unit 
${\1} = 0$. Denote this semiring by $\Rmax$; it is idempotent. 
The semiring $\Rmax$ is actually a 
commutative semifield. This construction is due to V.P. Maslov \cite{1};
now it is known as {\it Maslov's dequantization}.

The analogy with quantization is obvious; the parameter $h$ plays the
r\^{o}le of the Planck constant, so $\Rpl$ (or $\bR$) can be viewed as
a ``quantum object'' and $\Rmax$ as the result of its ``dequantization''. A
similar procedure gives the semiring ${\Rmin} = {\bR} \cup \{+\infty\}$ with
the operations $\oplus = \min$, $\odot = +$; in this case ${\0}=+\infty$, 
${\1} = 0$. The semirings $\Rmax$ and $\Rmin$ are isomorphic. Connections
with physics and imaginary values of the Planck constant are discussed
in \cite{4}. The commutative idempotent semiring ${\bR} \cup \{-\infty\} \cup
\{+\infty\}$ with the operations $\oplus = \max$, $\odot = \min$ can be
obtained as a result of a ``second dequantization'' of $\bR$ (or
$\Rpl$).  Dozens of interesting examples of nonisomorphic idempotent
semirings may be cited as well as a number of standard methods of deriving
new semirings from these (see, e.g., \cite{2}--\cite{5}).

{\it Idempotent dequantization} is a generalization of Maslov's dequantization. This
is a passage from fields to idempotent semifields and semirings in
mathematical constructions and results. The idempotent correspondence
principle (see Introduction and \cite{2, 4}) often works for this idempotent 
dequantization.

{\bf 3. Idempotent Analysis.}
Let $S$ be an arbitrary semiring with idempotent addition $\oplus$ (which
is always assumed to be commutative), multiplication $\odot$, zero $\0$, and
unit $\1$. The set $S$ is supplied with the {\it standard partial
order\/}~$\cle$: by definition, $a \cle b$ if and only if $a \oplus b = b$.
Thus all elements of $S$ are positive: ${\0} \cle a$ for all $a \in S$. Due
to the existence of this order, Idempotent Analysis is closely related to
the lattice theory, the theory of vector lattices, and the theory
of ordered spaces. Moreover, this partial order allows to
simulate a number of basic notions and results of Idempotent Analysis at
the purely algebraic level.

Calculus deals mainly with functions whose values are numbers. The
idempotent analog of a numerical function is a map $X \to S$, where $X$ is
an arbitrary set and $S$ is an idempotent semiring. Functions with values
in $S$ can be added, multiplied by each other, and multiplied by elements
of~$S$.

The idempotent analog of a linear functional space is a set of $S$-valued
functions that is closed under addition of functions and multiplication of
functions by elements of $S$, or an $S$-semimodule. Consider, e.g., the
$S$-semimodule ${\CB}(X, S)$ of functions $X \to S$ that are bounded in
the sense of the standard order on~$S$.

If $S = \Rmax$, then the idempotent analog of integration is defined by the
formula
\begin{equation}
   I(\varphi) = \int_X^{\oplus} \varphi (x)\, dx
	= \sup_{x\in X} \varphi (x),
\end{equation}
where $\varphi \in {\CB}(X, S)$. Indeed, a Riemann sum of the form
$\sumlim_i \varphi(x_i) \cdot \sigma_i$ corresponds to the expression
$\bigoplus\limits_i \varphi(x_i) \odot \sigma_i = \maxlim_i \{\varphi(x_i)
+ \sigma_i\}$, which tends to the right-hand side of~(1) as $\sigma_i \to
0$. Of course, this is a purely heuristic argument.
Formula~(1) defines the idempotent integral not only for functions taking
values in $\Rmax$, but also in the general case when any of bounded
(from above) subsets of~$S$ has the least upper bound.

An idempotent measure on $X$ is defined by $m_{\psi}(Y) = \suplim_{x \in Y}
\psi(x)$, where $\psi \in {\CB}(X,S)$. The integral with respect to this
measure is defined by
\begin{equation}
   I_{\psi}(\varphi)
	= \int^{\oplus}_X \varphi(x)\, dm_{\psi}
	= \int_X^{\oplus} \varphi(x) \odot \psi(x)\, dx
	= \sup_{x\in X} (\varphi (x) \odot \psi(x)).
\end{equation}

Obviously, if $S = \Rmin$, then the standard order $\cle$ is opposite to
the conventional order $\leqslant$, so in this case equation~(2) assumes the form
\begin{equation}
   \int^{\oplus}_X \varphi(x)\, dm_{\psi}
	= \int_X^{\oplus} \varphi(x) \odot \psi(x)\, dx
	= \inf_{x\in X} (\varphi (x) \odot \psi(x)),
\end{equation}
where $\inf$ is understood in the sense of the conventional order $\leqslant$.

The functionals $I(\varphi)$ and $I_{\psi}(\varphi)$ are linear over $S$;
their values correspond to limits of Lebesgue (or Riemann) sums.
The formula for $I_\psi(\varphi)$ defines the
idempotent scalar product of the functions $\psi$ and $\varphi$. Various
idempotent functional spaces and an idempotent version of the theory of
distributions can be constructed on the basis of the idempotent integration,
see, e.g., \cite{1}, \cite{3}--\cite{5}. The analogy of idempotent and probability measures 
leads to spectacular parallels between optimization theory and probability 
theory. For example, the Chapman--Kolmogorov equation corresponds to the 
Bellman equation (see, e.g., \cite{6, 5}). Many other idempotent analogs may 
be cited (in particular, for basic constructions and theorems of functional 
analysis~\cite{4}).

{\bf 4. The superposition principle and linear problems.}
 Basic equations of quantum theory are linear (the superposition principle). 
The Hamilton--Jacobi equation, the basic equation of classical mechanics, is 
nonlinear in the conventional sense. However it is linear over the semirings 
$\Rmin$ and $\Rmax$. Also, different versions of the Bellman equation, the 
basic equation of optimization theory, are linear over suitable idempotent
semirings (V.~P.~Maslov's idempotent superposition principle), see \cite{1, 3}.
For instance, the finite-dimensional stationary Bellman equation can be 
written in the form $X = H \odot X \oplus F$, where $X$, $H$, $F$ are 
matrices with coefficients in an idempotent semiring $S$ and the unknown 
matrix $X$ is determined by $H$ and $F$ \cite{7}. In particular, standard 
problems of dynamic programming and the well-known shortest path problem 
correspond to the cases $S = {\Rmax}$ and $S ={\Rmin}$, respectively. 
In \cite{7}, it was shown that main optimization algorithms for finite graphs 
correspond to standard methods for solving systems of linear equations of this 
type (i.e., over semirings). Specifically, Bellman's shortest path algorithm
corresponds to a version of Jacobi's algorithm, Ford's algorithm
corresponds to the Gauss--Seidel iterative scheme,~etc.

Linearity of the Hamilton--Jacobi equation over $\Rmin$ (and $\Rmax$) is
closely related to the (conventional) linearity of the Schr\"odinger equation,
see \cite{4} for details.

{\bf 5. Correspondence principle for algorithms and their computer
implementations.} The idempotent correspondence principle is valid
for algorithms as well as for their software and hardware implementations
\cite{2}. In particular, according to the superposition principle, analogs
of linear algebra algorithms are especially important. It is well-known
that algorithms of linear algebra are convenient for parallel computations;
so their idempotent analogs accept a parallelization. This is a regular way 
to use parallel computations for many problems including basic optimization
problems. It is convenient to use universal algorithms which do not
depend on a concrete semiring and its concrete computer model. Software
implementations for universal semiring algorithms are based on 
object-oriented and generic programming; program modules can deal with 
abstract (and variable) operations and data types, see \cite{2, 8} for
details.

The most important and standard algorithms have many hardware 
realizations in the form
of technical devices or special processors. These devices often can be used as 
prototypes for new hardware units generated by substitution of the
usual arithmetic operations for its semiring analogs, see \cite{2} for
details. Good and efficient technical ideas and decisions can be 
transposed from prototypes into new hardware units. Thus the correspondence
principle generates a regular heuristic method for hardware design. Note 
that to get a patent it is necessary to present the so-called 
``invention formula'', that is to indicate a prototype for the suggested
device and the difference between these devices.

{\bf 6. Idempotent interval analysis.} An idempotent version of the
traditional interval analysis is presented in \cite{9}. Let $S$ be an
idempotent semiring equipped with the standard partial order (see the
beginning of Section 3). A {\it closed interval} in $S$ is a subset of the 
form ${\bx} = [{\bxu}, {\bxo}] = \{ x\in S\vert {\bxu}\cle x \cle{\bxo}\}$, 
where the elements $\bxu\cle \bxo$ are called {\it lower} and {\it upper 
bounds} of the interval $\bx$. A {\it weak interval extension} $I(S)$ of the
semiring $S$ is the set of all closed intervals in $S$ endowed with
operations $\oplus$ and $\odot$ defined as $\bx\oplus\by = [\bxu\oplus\byu, 
\bxo\oplus\byo]$, $\bx\odot\by = [\bxu\odot\byu, \bxo\odot\byo]$; the set
$I(S)$ is a new idempotent semiring with respect to these operations.
It is proved that basic problems of idempotent linear algebra are
polynomial, whereas in the traditional interval analysis problems of
this kind are generally NP-hard. Exact interval solutions for the
discrete stationary Bellman equation (this is the matrix equation
discussed in Section 4) and for the corresponding optimization 
problems are constructed and examined.

{\bf 7. Generalized fuzzy sets.} Let $\Omega$ be the so-called universe
consisting of ``elementary events'' and $S$ an idempotent semiring.
Denote by ${\CF}(S)$ the set of functions defined on $\Omega$ and taking
their values in $S$; then ${\CF}(S)$ is an idempotent semiring with respect
to the pointwise addition and multiplication of functions. We shall
say that elements of ${\CF}(S)$ are {\it generalized fuzzy sets}. See also \cite{13}. We have the 
well known classical definition of fuzzy sets (L.A. Zadeh \cite{10})
if $S = {\bP}$, where $\bP$ is the segment $[0,1]$ with the semiring operations
$\oplus = \max$ and $\odot = \min$, see Section 2. Of course,
functions from ${\CF}(\bP)$ taking their values in the Boolean algebra
${\bB} = \{0, 1\}\subset {\bP}$ correspond to traditional sets from $\Omega$
and semiring operations correspond to standard operations for sets. In the
general case if $S$ has neutral elements $\0$ and $\1$ (and ${\0}\neq {\1}$),
then functions from ${\CF}(S)$ taking their values in $\bB = \{{\0}, 
{\1}\}\subset S$ can be treated as traditional subsets in $\Omega$. If $S$
is a lattice (i.e. $x\odot y = \inf \{x, y\}$ and $x\oplus y = \sup 
\{x, y\}$), then generalized fuzzy sets coincide with $L$-fuzzy sets in the
sense of J.A. Goguen \cite{11}. The set $I(S)$ of intervals is an idempotent
semiring (see Section 6), so elements of ${\CF}(I(S))$ can be treated as
interval (generalized) fuzzy~sets.

It is well known that the classical theory of fuzzy sets is a basis
for the theory of possibility \cite{12}. Of course, it is possible
to develop a similar generalized theory of possibility starting from
generalized fuzzy sets. In general the generalized theories are
noncommutative; they seem to be more qualitative and less quantitative
with respect to the classical theories presented in \cite{10, 12}. We see that
Idempotent Analysis and the theory of (generalized) fuzzy sets and 
possibility have the same objects, i.e. functions taking their values in
semirings. However, basic problems and methods could be different
for these theories (like for the measure theory and the probability
theory).

\end{document}